\documentclass[preprint,sort,compress,11pt]{elsarticle}

\usepackage{amssymb}
\usepackage{amsthm}
\usepackage{amsmath}
\usepackage{booktabs}
\usepackage{multirow}

\usepackage{url}
\usepackage[all]{nowidow} % tex.stackexchange.com/a/30409
\usepackage[colorlinks=true, linkcolor=blue, urlcolor=blue, citecolor=blue]{hyperref}

\allowdisplaybreaks

\DeclareMathOperator*{\argmax}{arg\,max}

\journal{Omega}

\begin{document}

\begin{frontmatter}



\title{Integrated optimization of operations and capacity planning under uncertainty for drayage procurement in container logistics}


\author[inst1,inst2]{Georgios Vassos}

\affiliation[inst1]{organization={Transported by Maersk, A.P. Moller - Maersk},%Department and Organization
            addressline={Esplanaden 50}, 
            city={Copenhagen~K},
            postcode={1098}, 
            state={},
            country={Denmark}}

\author[inst2]{Richard Lusby}

\affiliation[inst2]{organization={Department of Technology, Management and Economics, Technical University of Denmark},%Department and Organization
            addressline={Akademivej, 358}, 
            city={Kgs.~Lyngby},
            postcode={2800}, 
            state={},
            country={Denmark}}

\author[inst3,inst2]{Pierre Pinson}

\affiliation[inst3]{organization={Dyson School of Design Engineering, Imperial College London},%Department and Organization
            addressline={Dyson Building, 1M04A}, 
            city={London},
            postcode={SW7~2AZ}, 
            state={},
            country={United Kingdom}}

% \affiliation[*]{Corresponding author}

\begin{abstract}
    We present an integrated framework for truckload procurement in container logistics, bridging strategic and operational aspects that are often treated independently in existing research. Drayage, the short-haul trucking of containers, plays a critical role in intermodal container logistics. Using dynamic programming, we identify optimal operational policies for allocating drayage volumes among capacitated carriers under uncertain container flows and spot rates. The computational complexity of optimization under uncertainty is mitigated through sample average approximation. These optimal policies serve as the basis for evaluating specific capacity arrangements. To optimize capacity reservations with strategic and spot carriers, we employ an efficient quasi-Newton method. Numerical experiments demonstrate significant cost-efficiency improvements, including a 21.2\% cost reduction in a four-period scenario. Monte Carlo simulations further highlight the strong generalization capabilities of the proposed joint optimization method across out-of-sample scenarios. These findings underscore the importance of integrating strategic and operational decisions to enhance cost efficiency in truckload procurement under uncertainty.
\end{abstract}


% \begin{highlights}
%     \item Novel framework integrates capacity planning with volume allocation policies in container drayage procurement.
%     \item Multidimensional uncertainties: inflow, outflow, and spot market variability.
%     \item Approximation methods mitigate computational complexity of dynamic programming.
%     \item Model generalizes robustly across unseen scenarios.
% \end{highlights}


\end{frontmatter}


\section{Introduction}

We study the complexities of \textit{truckload (TL)} procurement for the short-haul transportation of laden containers, a critical process in managing the first-mile movement of \textit{Full Container Load (FCL)} and \textit{Less-than-Container Load (LCL)} shipments between ports, warehouses, or intermodal terminals. This stage, often referred to as \textit{container drayage} or simply \textit{drayage}, plays a pivotal role in bridging ocean freight with inland logistics \citep{Nossack2013}. Throughout this work, we use the terms TL and drayage interchangeably, as TL, while a broader category, refers exclusively to drayage in this context.

\textit{First mile delivery} refers to the initial phase of the logistics chain, where goods are collected from the shipper’s location and transported to a warehouse, sorting facility, or terminal as the starting point for further distribution. It primarily concerns the movement of goods from the point of origin to the beginning of the transportation network. FCL is a shipping mode where a container is exclusively used for one customer's cargo, filling the entire container. In contrast, LCL is another shipping arrangement where cargo from multiple customers is consolidated into a single container. \textit{Intermodal container logistics} refers to the integrated transportation of containers using multiple modes of transportation, such as vessels, trucks, and trains, without any handling of the freight itself when changing modes. Finally, drayage involves the short-distance trucking of containerized freight into and out of intermodal facilities, such as ports, rail terminals, and warehouses, facilitating the transition between transportation modes.

Efficient TL procurement for drayage is crucial for maintaining smooth intermodal logistics operations, minimizing costs, and enhancing service levels. As global trade grows in complexity and volume, optimizing these processes has substantial implications for supply chain efficiency and economic performance \citep{Caplice2006c}. However, much of the existing research focuses on strategic carrier selection and auction mechanisms, often overlooking the critical role of capacity planning \citep{acocella2023a}.

In this study, we distinguish between carrier selection and capacity planning as two distinct yet complementary components of the procurement process. Carrier selection entails identifying and contracting carriers based on criteria such as cost, reliability, and service quality \citep{Taherdoost2019}, establishing the strategic foundation for aligning carriers with operational needs. In contrast, capacity planning focuses on determining the transportation capacity to reserve with both strategic and spot carriers, ensuring adequate resources are available to meet anticipated shipment volumes efficiently.

While carrier selection has received significant attention for fostering long-term partnerships in TL transportation, capacity planning remains underexplored. Neglecting this component can impede the integration of strategic and operational considerations in TL procurement. Effective procurement should secure sufficient capacity with both strategic and spot carriers to meet operational needs \citep{Martinez2005}. This gap in the literature underscores the need for integrated models that can enhance strategic decision-making in TL procurement.

Our study integrates strategic capacity planning, which involves the long-term reservation of carrier capacity to accommodate expected container flows, with operational volume allocation, which entails the periodic management of transportation volumes. Given a predefined capacity arrangement, we aim to identify the optimal policy for periodically assigning transportation volumes across capacitated carriers under uncertain flows of laden containers within the drayage system. To achieve this, we combine the well-established capacitated transportation problem \citep{Castillo2001, Korman2015} with dynamic programming. Additionally, we explore approximate methods to manage computational complexity, emphasizing its strong dependence on uncertain parameters. % Finally, we specify the regularity conditions required to ensure the capacity search problem is concave, enabling it to be efficiently solved using first-order gradient methods.

To address the gap in literature, we propose a joint optimization model, combining capacity planning with operational policies for volume allocation. We use the \textit{Markov Decision Process (MDP)} framework to formulate a dynamic program for the \textit{Multiperiod Stochastic Transportation Problem (MSTP)}. This method can readily factor in multiple uncertainties, such as the inflow and outflow of laden containers to and from our drayage system and the spot market variability. Our approach is convenient for practitioners who wish to experiment with empirical data, as the uncertain parameters can be fine-tuned and validated to match real-world mechanisms. 

Additionally, we adopt the \textit{option contract} framework and the \textit{portfolio contract} concept introduced by \cite{Martinez2005} to enable flexible sourcing, extending the capacity strategy to include the spot market. An option contract grants the buyer the right, but not the obligation, to procure a specified amount of capacity at a predetermined execution price in exchange for an upfront reservation fee. A portfolio contract combines multiple option contracts, providing flexibility to optimize operational efficiency. Finally, we leverage the theory therein to report the regularity conditions necessary to reduce the capacity optimization at the outset of operations to a concave maximization problem. 

The remainder of this paper is structured as follows: Section~\ref{sec:litrev} reviews key contributions to truckload procurement and identifies existing gaps. Section~\ref{sec:probdescr} defines essential concepts and presents the mathematical formulation of the optimization model. Section~\ref{sec:method} introduces an MDP framework for optimizing volume allocation under uncertainty, integrating strategic capacity planning. Computational complexity and scalability are analyzed using dynamic programming and approximation techniques. Section~\ref{sec:numex} presents numerical experiments that demonstrate the model's effectiveness, followed by a discussion of the results, limitations, and future research directions.

% \vspace*{\fill}

\section{Literature review}\label{sec:litrev}

The procurement of TL transportation services has been extensively studied from the perspectives of strategic carrier selection and operational execution. This section reviews the state-of-the-art research, highlighting key contributions and identifying gaps that this study seeks to address.

Strategic carrier selection in TL procurement involves making long-term decisions about which carriers to contract and the terms of these contracts. This topic has garnered significant attention, with numerous studies exploring different aspects of the process. Researchers have developed models to determine the optimal set of carriers and allocate them to \textit{origin-destination (OD)} pairs, highlighting the importance of auction mechanisms in achieving cost efficiency and maintaining service quality \citep{Caplice2003}. Further investigations into the TL reverse auction process have provided insights into bid construction and the complexities of carrier selection, showing how bid strategies influence procurement auction outcomes \citep{Caplice2006c, Caplice2007e}. Subsequent studies have employed multi-attribute value theory and auction theory to address scenarios such as lane bundling and uncertain demand and capacity \citep{Buer2014, Hu2016, Lim2008, Ma2010, Oner2021, Rekik2012, Remli2013, Remli2019}. In this context, a \textit{lane} refers to a specific route connecting an OD. Additionally, contributions to central auctioneer solutions and carrier assignment algorithms have expanded the understanding of carrier selection \citep{Sandholm2002, Xu2013, Xu2014, Zhang2014, Zhang2015}. Reviews of strategic TL procurement literature have also highlighted simplifying assumptions about demand and the limited focus on non-price factors, such as environmental sustainability \citep{Basu2015, Basu2017}.

Strategic management of carrier capacity is crucial for optimizing TL transportation procurement, ensuring operational efficiency and a competitive advantage. The literature has explored the optimization of carrier assignments while accounting for capacity constraints \citep{Caplice2003}. Algorithms have been developed to manage carrier capacity under contractual agreements in conditions of uncertain demand \citep{Lim2008, Zhang2015}. Risk mitigation strategies that align carrier capacity with fluctuating demand to maintain consistent service levels have also been proposed \citep{Yoon2016}. Furthermore, pricing strategies that integrate carrier capacity into bid generation for lane bundles have been examined \citep{MesaArango2015}.

Operational execution focuses on managing transportation services, considering contract commitments, carrier capacity, and uncertainty in shipment demands and market conditions. Research has emphasized the critical impact of market conditions and competitive contract prices on carrier acceptance \citep{Acocella2020}. Studies have demonstrated the benefits of collaborative strategies, such as reducing deadhead miles and enhancing negotiating power \citep{Ergun2007, Ozener2008}. Trade-offs between spot sourcing and strategic delays to maximize revenues have also been explored, indicating that while relying on the spot market can improve fill rates, it may compromise reliability and resource utilization \citep{Kantari2021, Mes2009}. The potential of option contracts to provide flexibility and reduce risks in TL transportation has been highlighted \citep{TibbenLembke2006}. The dynamic nature of TL operations, where partial demand information and market conditions influence contract commitments and capacity availability, underscores the importance of flexibility in decision-making \citep{Boada2020a}. Despite these advances, the influence of market uncertainty on procurement strategies remains underexplored. While some recommendations advocate for accommodating demand surges and cyclical market conditions, practical models to address these factors are still lacking \citep{acocella2023a, Vos1999}.

Despite extensive research on strategic carrier selection and operational execution, a significant gap remains in integrating these two phases \citep{acocella2023a}. Most studies address them independently, leading to inefficiencies and missed optimization opportunities. The need for integrated capacity planning is highlighted in studies such as \citep{Acocella2020, Acocella2022a, Acocella2022b}, which examine market conditions and contract price competitiveness. These works suggest that more dynamic and adaptive models are necessary to bridge the gap between strategic planning and operational execution. \cite{boujemaa2022sddp} introduced a multistage stochastic optimization model for carrier evaluation under uncertain demand in middle-mile logistics, focusing on TL transportation between warehouses and distribution centers. Our work extends this framework by incorporating uncertainties in inflow, outflow, and spot rates, as well as a strategic capacity planning optimization stage.

To sum up, the procurement of TL transportation services has been analyzed through the lenses of strategic carrier selection and operational execution. However, effectively integrating these critical phases remains a challenge. The coupling of carrier selection and capacity planning adds further complexity to this integration. Moreover, there is a notable absence of holistic models capable of adapting to market conditions and the uncertainties of supply and demand. To address this gap, this study proposes a framework for drayage procurement that incorporates operational impacts into strategic decision-making under uncertainty in container flows and spot market.

\section{Problem description}\label{sec:probdescr}

We describe the transportation problem within the drayage process in the context of intermodal container logistics. This section introduces the key categories required for understanding this problem as a practitioner in the container logistics industry. It also highlights the role of TL procurement in the drayage process and underscores the importance of an effective capacity strategy. 

Intermodal container logistics is a complex process that involves efficiently moving containers using multiple transportation modes, such as ships, trains, and trucks. This process encompasses the first mile, which involves transporting containers to an intermodal terminal; FCL, where fully loaded containers are booked by individual customers; and LCL, where shipments from multiple customers are consolidated into a single container for the primary transport leg. Drayage is critical for facilitating modality transitions, employing short-distance trucking to transport containers from arrival points to nearby destination hubs.

\begin{figure}[ht]
    \centering
    \includegraphics[width = 0.75 \textwidth]{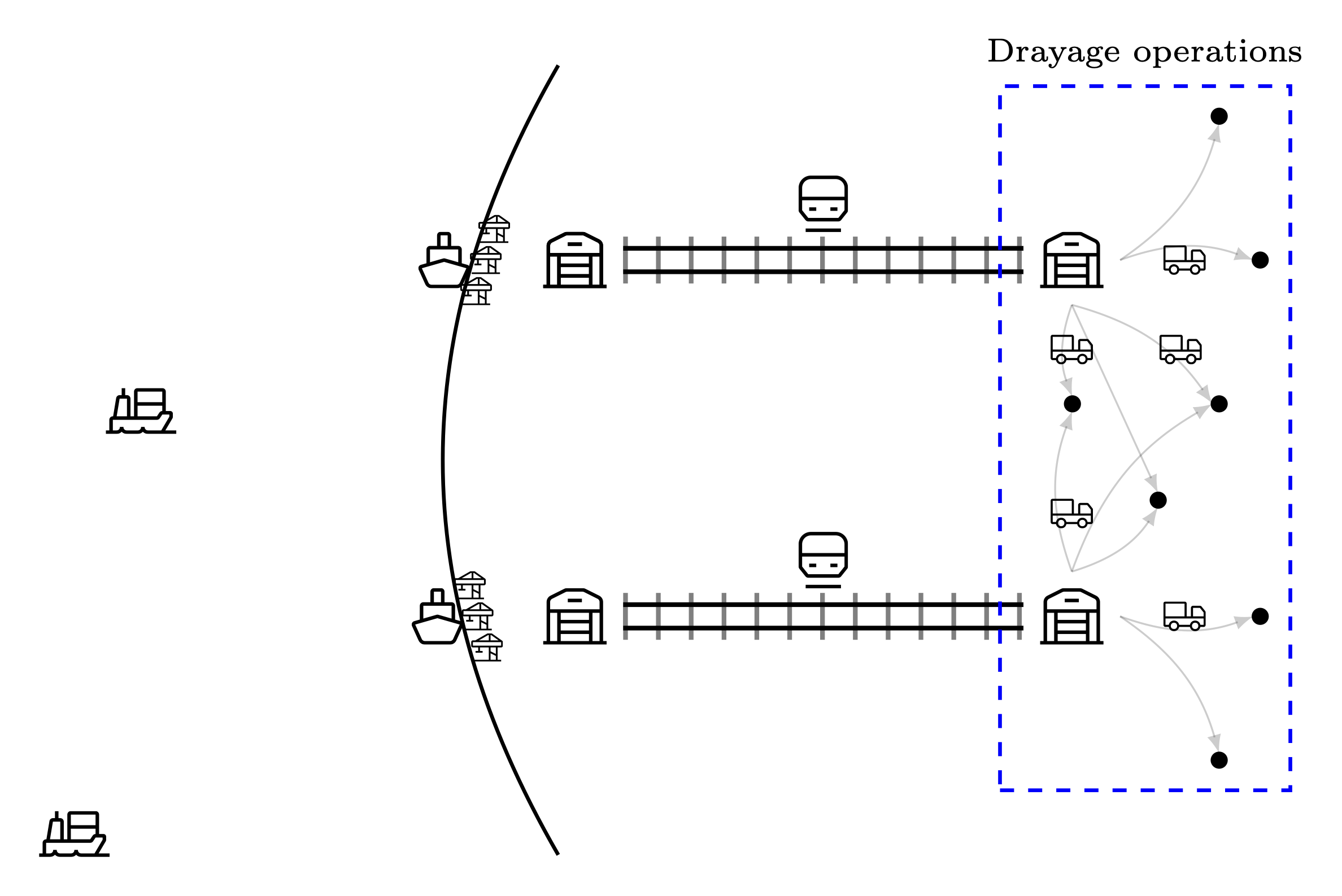}
    \caption{An example of intermodal container transportation involves containers being transported by ocean and rail into a system of drayage operations. These containers arrive at rail yards within the system's scope, and from there, they are transported to nearby hubs, where they continue their intermodal journey.}
    \label{fig:drayops}
\end{figure}

A system of drayage operations includes storage locations at entry points, such as ports or rail yards, where containers arrive following an ocean or rail leg. These containers are subsequently transported by trucks to neighboring exit points, which may include other ports, rail yards, or intermediate hubs. Figure~\ref{fig:drayops} illustrates a typical drayage operation within container logistics. Containers arrive at entry points from a combined ocean-rail leg and are transported by truck to adjacent hubs to continue their journey. This system is characterized by a continuous inflow and outflow of containers, which are exogenous to the drayage system and beyond its direct control.

In this study, we examine a case where a shipping company outsources its entire TL transportation business to external carriers. This is managed through a procurement process that begins with the onboarding of strategic carriers, referred to as carrier selection \citep{acocella2023a}, and includes auctions where carriers compete for bundles of lanes \citep{Buer2014, MesaArango2015, Oner2021}, also known as corridors in intermodal container logistics. The outcome of carrier selection establishes the foundation of the procurement strategy through long-term service agreements with multiple carriers. During this process, carriers bid for specific corridors, with winning bids specifying fixed trucking rates per kilometer.

Capacity planning involves determining the trucking capacity to secure with strategic carriers and setting buffer sizes for the spot market, where transportation services are procured to meet immediate needs in a competitive environment. Buffers ensure adaptability to container flow volatility by reserving capacity. To enable flexible sourcing, we incorporate option contracts and the portfolio contract introduced by \cite{Martinez2005}. Option contracts involve pre-paying a premium to reserve capacity, granting the right to procure up to a specified amount at a fixed price. The portfolio contract builds on this by integrating long-term, option, and spot market agreements, optimizing procurement to balance cost and flexibility under uncertainty.

Overall, this study addresses periodic transportation decisions related to the movement of container volumes from entry to exit points, the allocation of business among capacitated carriers through the portfolio contract, and the capacity strategy that sets carrier capacity limits for the planning horizon. To provide a clearer understanding of these elements, the next section describes an elementary system of drayage operations and presents the mathematical notation conventions used throughout this work.

\subsection{Mathematical formulation and notation}

Let \(\mathcal{I}\) and \(\mathcal{J}\) be the index sets for entry and exit points, respectively, and let \( \mathcal{T} \) represent the index set for the planning horizon, spanning from time point \( 1 \) until the stopping time \( \tau < \infty \). We define \(Q_t = (Q_{i,t} : i \in \mathcal{I})\) as the inflow of containers entering our system at time \(t\), and \(D_t = (D_{j,t} : j \in \mathcal{J})\) as the outflow of containers exiting. Transportation within our system occurs through lanes denoted by \(\mathcal{L} \subset \mathcal{I} \times \mathcal{J}\). The state vector \(S_t = (S_{i,t} : i \in \mathcal{I}, S_{j,t} : j \in \mathcal{J})\) represents the number of containers stored at each entry and exit location in our system at time \(t\). The decision vector \( \check{A}_t = (A_{i,j,t} : (i, j) \in \mathcal{L}) \) represents the total volume of containers transported through each lane \((i, j)\) at time \(t\).

\begin{figure}[ht]
    \centering
    \includegraphics[width = 0.5\textwidth]{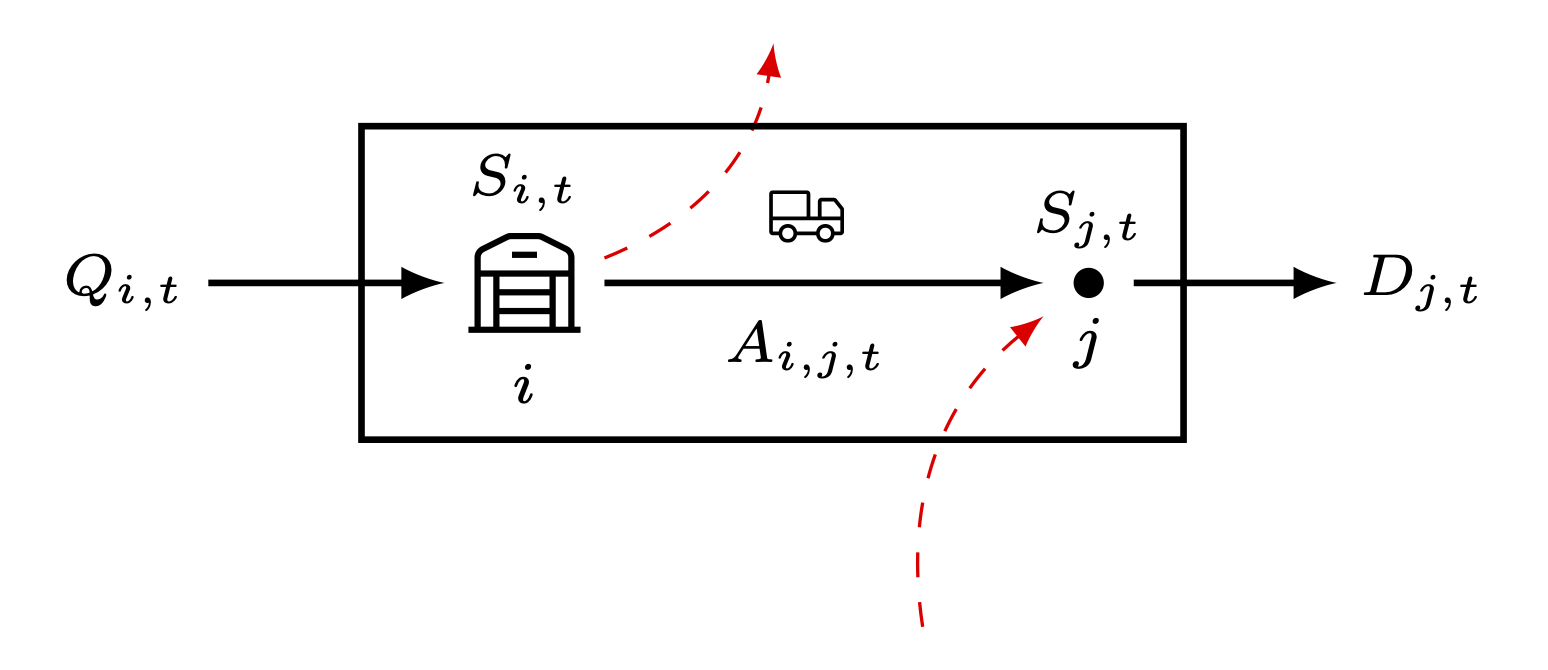}
    \caption{An elementary drayage system with storage facilities at the entry and exit points.}
    \label{fig:elementary}
\end{figure}

Figure~\ref{fig:elementary} illustrates the notation using an elementary snippet of a drayage system that comprises a single lane. At time \( t \), the state of the system is described by the number of containers stored at the entry and exit hubs, \( (S_{i,t}, S_{j,t}) \). There is an inflow of \( Q_{i,t} \) containers and an outflow of \( D_{j,t} \) containers. Additionally, we may intervene in the system by transporting \( A_{i,j,t} \) containers from \( i \) to \( j \).

\subsubsection{System dynamics}

In a cycle of operations spanning a finite horizon indexed by $t \in \mathcal{T}$, where $\mathcal{T} = \{1, 2, \dots, \tau\}$, our system evolves according to the following sequential logic. At each time stage $t$, the drayage system is at state $S_{t}$, quantities $\sum_{i\in\mathcal{I}}Q_{i,t}$ of containers arrive to entry points, the decision-maker orders a number of drayage moves $A_{t}$, and a volume of $\sum_{j\in\mathcal{J}}D_{j,t}$ containers departs from exit points. What lies latent in the realizations of the exogenous variables is the dynamics of the broader container flow in the company's network of intermodal container transportation.

We now introduce several concepts that will allow us to formalize the representation of the uncertainty in our problem. Let \((\Omega, \mathcal{F}, P)\) be a probability space, where \(\Omega\) is the sample space of all possible realizations of the uncertainty, \(\mathcal{F}\) is the \(\sigma\)-algebra containing all measurable subsets of \(\Omega\), and \(P: \mathcal{F} \to [0, 1]\) is a probability measure assigning a probability to each measurable set. In this context, we may refer to a measurable set as an event. A random variable \(Z\) is a measurable function \(Z: \Omega \to \mathcal{Z}\) defined on our probability space, where \( \mathcal{Z} \) is the space of all possible outcomes of \( Z \). 

Building on this foundation, a stochastic process is a collection of random variables indexed by time, denoted \((Z_t : t \ge 0)\). To capture the progression of information, we define a filtration as an increasing sequence of sub-\(\sigma\)-algebras of \(\mathcal{F}\), indexed by time, denoted by \((\mathcal{F}_t : t \ge 0)\). A stochastic process \((Z_t : t \ge 0)\) is adapted to a filtration \((\mathcal{F}_t : t \ge 0)\) if \(Z_t\) is \(\mathcal{F}_t\)-measurable for every \(t \ge 0\). This provides a minimal framework to represent and analyze uncertainty in our problem.

In this work, we consider the discrete-time setting and define $Z_t : \Omega \to \mathcal{Z}$ to represent the uncertainty at time $t$. For now, $Z_t = (Q_t, D_t)$, where $Q_t$ and $D_t$ denote the random inflows and outflows, respectively; however, $Z_t$ will later be extended to incorporate spot rates. To formalize the exogenous state, we recursively define $\Phi_t = (\Phi_{t-1}, Z_t)$ for all $t \geq 1$, with $\Phi_0 = \emptyset$. Here, $\Phi_t$ captures the history of the uncertainty up to time $t$, and is a random variable $\Phi_t : \Omega \to \mathcal{Z}^t$, where the set $\mathcal{Z}^t$ is defined as the $t$-fold Cartesian product of $\mathcal{Z}$:
\[
\mathcal{Z}^t = \mathcal{Z} \times \mathcal{Z} \times \cdots \times \mathcal{Z} \quad (t \text{ times}),
\]
with $\mathcal{Z}^0 = \{\emptyset\}$. Thus, $\mathcal{Z}^t$ represents the set of all possible histories up to time $t$.

We fix a history $\phi_{\tau} \in \mathcal{Z}^{\tau}$, and the state of our system evolves conditional on this fixed realization. For brevity, we omit the conditioning on the fixed history in our transition probability and let $p_{t}(s_{t+1} \mid s_{t}, a_{t}) \doteq p_{t}(s_{t+1} \mid s_{t}, a_{t}, \phi_{\tau})$. After an action \( a_{t} \in \mathcal{A}(s_{t}) \) is carried out in a given state \( s_{t} \), the system transitions to a new state according to the deterministic law:
\begin{equation}\label{eq:draydyn}
    p_{t}(s_{t+1} \mid s_{t},a_{t}) = \delta(s_{t+1}-f_{S,t}(s_{t},a_{t},\phi_{t}))
\end{equation}
Here, \( \delta \) denotes Dirac's delta, a generalized function that assigns nonzero weight only when its argument is zero, formalizing the deterministic nature of the transition \citep{bracewell1999}. The function $f_{S,t} : \mathcal{S} \times \mathcal{A} \times \mathcal{Z}^{t} \to \mathcal{S}$ in the structural equation $S_{t+1}=f_{S,t}(S_{t},A_{t},\Phi_{t})$, which determines the next state, encapsulates the system dynamics by providing the mechanism through which the system will respond over time. It is defined as follows:
\begin{equation}
    f_{S,t}(S_{t},A_{t},\Phi_{t})\doteq\begin{Bmatrix}
        \left(S_{i,t}-\sum_{j\in\mathcal{J}}A_{i,j,t}+Q_{i,t}:i\in\mathcal{I}\right),\\[5pt]
        \left(S_{j,t}+\sum_{i\in\mathcal{I}}A_{i,j,t}-D_{j,t}:j\in\mathcal{J}\right)
    \end{Bmatrix}
\end{equation}
For a given instance $\phi_{t} \in \mathcal{Z}^{t}$ of the uncertainty, the transition to the next state is deterministic. The only nuisance parameter in the transition mechanism is the density function of the exogenous state $p_{t}(\phi_{t}) = P_{t}(\Phi_{t}=\phi_{t})$. The complexity of the exogenous state variable combined with the lack of sufficient data, poses a great challenge in validating any modeling choice for this nuisance parameter. Nevertheless, we can still gain valuable insights into drayage operations by simulating the nuisance parameter and analyzing the evolution of our system.

\subsubsection{Flexible sourcing tools}

During the carrier selection process, a set of bids $\mathcal{B}$ is considered, each associated with a subset of lanes $\mathcal{L}(b) \subset \mathcal{L}$, a fixed transportation cost, and potentially, a one-time reservation cost. The union of these subsets covers all available lanes, i.e., $\mathcal{L}(\mathcal{B}) = \bigcup_{b \in \mathcal{B}} \mathcal{L}(b) = \mathcal{L}$. Carriers with winning bids are then designated as strategic carriers, denoted by $\mathit{CS}$. For each strategic carrier, identified as $k \in \mathit{CS}$, $\mathcal{B}(k) \subset \mathcal{B}$ represents the set of bids that have been successfully awarded to them. Furthermore, $\mathcal{L}(k) = \bigcup_{b \in \mathcal{B}(k)} \mathcal{L}(b)$ defines the specific lanes that each carrier is responsible for managing. To simplify the notation, we mildly abuse the term \textit{source} to refer to the bid and winner carrier pair \((b,k)\), and use a single index $k$ to identify a source. This is necessary because a carrier may win the same lanes under different bids, each with distinct terms, and treating the (bid, carrier) pair as the unit ensures we capture those differences while simplifying the notation.

A typical long-term service agreement involves a carrier's commitment to a fixed transportation cost for each lane within the bid, coupled with a capacity assurance. Such an agreement for a given source $k$ can be represented as the vector $\left((w_{t}^{k}, x_{t}^{k}) : t \in \mathcal{T}\right)$, where $w_{t}^{k}$, defined as $\left( w_{i,j,t}^{k} : (i, j) \in \mathcal{L}(k) \right)$, represents the execution cost at time $t$ for each lane under the purview of source $k$. The variable $x_{t}^{k}$, taking values in $\mathbb{R}^{+} = [0, \infty)$, denotes the capacity reserved with source $k$ at time $t$. The variation in transportation costs among lanes is justifiable due to differences in distance coverage and the involvement of distinct bids.

An option contract is a flexible sourcing tool used for managing uncertainties in supply chain operations. For a given source, it is defined by the vector $\left( (v_{t}^{k}(\,\cdot\,),w_{t}^{k},x_{t}^{k}):t\in\mathcal{T} \right)$, where $v_{t}^{k}:\mathbb{R}^{+}\to\mathbb{R}^{+}$ is the premium cost function for reserving a specified amount of capacity with the source $k$ for time period $t$. The variables $w_{t}^{k}$ and $x_{t}^{k}$ represent the fixed execution cost and the reserved capacity, respectively, at time $t$, similar to the terms found in service agreements. Option contracts have been proposed to effectively mitigate both price and demand risks in TL procurement and other areas \citep{Anderson2017INF, Cai2016, Martinez2005, Tsai2011}. Moreover, integrating spot sources into this strategy is a straightforward extension. Let $\mathit{CO}$ denote the set of available spot carriers. A spot source $k\in\mathit{CO}$ can be modeled as an option contract with no premium cost, offering high capacity and random execution rates $(W_{t}^{k}:t\in\mathcal{T})$, which reflect the typical characteristics of spot market transactions.

% A portfolio contract in this setting consists of a structured combination of option contracts across multiple sources. This framework enables shippers to balance cost efficiency with capacity flexibility by selecting a mix of long-term service agreements and spot market options. By structuring the portfolio in this manner, firms can hedge against price fluctuations and demand uncertainty, ensuring both cost competitiveness and service reliability.

\subsection{Objective function}

The primary objective is to determine optimal volume allocation policies for a set of capacitated sources over the planning horizon. Let \( \pi = (\pi_{t} : t \in \mathcal{T}) \) denote a policy, where \( \pi_{t} : \mathcal{S} \times \mathcal{A} \to [0,1] \) defines a stochastic policy, and \( \pi_{t} : \mathcal{S} \times \mathcal{A} \to \{0, 1\} \) defines a deterministic policy for all \( t \in \mathcal{T} \). For every \( t \in \mathcal{T} \), let \( V_{t} : \mathcal{S} \to \mathbb{R} \) represent the value function, which assigns a numerical value to each system state at time \( t \). The value function under a policy \( \pi \in \Pi \), where \( \Pi \) is a set of policies, is denoted by \( V_{t}^{\pi} \). For any \( s_{1} \in \mathcal{S} \), the optimal policy \( \pi^{*}(s_{1}) = ( \pi_{1}(s_{1}), \pi_{t} : t = 2, \dots, \tau ) \), with \( \pi_{1}(s_{1}):\mathcal{A} \to \{0, 1\} \), is a deterministic policy defined as the policy that maximizes the value of the system's initial state, given by:
\begin{equation}
    \pi^{*}(s_{1}) = \argmax_{\pi \in \Pi} V_{1}^{\pi}(s_{1})
\end{equation}
This can be conditional of a specific scenario \( \phi_{\tau} = (z_{t} : t \in \mathcal{T}) \) or across all scenarios by optimizing the expectation functional \( E_{\Phi_{\tau}} V_{1}^{\pi}(s_{1}) \) instead. We provide more details on this in the next section, where the value function for the volume allocation problem in drayage is defined. We take into account uncertainties in the inflow and outflow of containers, as well as the variable costs associated with the spot market. 

The value of an initial state \( s_{1} \) under the optimal policy \( \pi^{*} \doteq \pi^{*}(s_{1}) \) is defined as the maximum achievable value at that state:
\begin{equation}
    V_{1}^{*}(s_{1}) = V_{1}^{\pi^{*}}(s_{1}) = \max_{\pi \in \Pi} V_{1}^{\pi}(s_{1})
\end{equation}
Ultimately, our goal is to determine the optimal capacities to reserve with the strategic carriers, represented by the index set \( \mathcal{K} \) of all carrier sources, and to establish the appropriate size of the spot market buffer, thereby ensuring effective drayage operations from the outset. Carrier capacities, denoted by \( x = (x_{t}^{k} : k \in \mathcal{K}, t \in \mathcal{T}) \), are incorporated as constraints in the value function computation. To reflect this dependency, we denote the value function under a specific capacity arrangement \( x \) as \( V_{1}^{*}(s_{1}; x) \).

For a fixed scenario \( \phi_{\tau} \in \mathcal{Z}^{\tau} \), the objective of the capacity optimization problem is to maximize \( V_{1}^{*}(s_{1}; x) - v(x) \), where \( v(x) \) represents the reservation cost of the capacity arrangement \( x \). In the following section, we develop a framework grounded in MDP theory to address these challenges.

\section{Method}\label{sec:method}

We propose an integrated optimization approach in which the first stage aims to optimize the capacity arrangement between strategic and spot carriers, while the second stage focuses on optimizing the operational policy for the drayage of container volumes from entry to exit points within our system.

\subsection{Optimal volume allocation policy}\label{sec:optops}

An MDP is a mathematical framework used for modeling decision-making scenarios where outcomes are influenced both by uncertainties and the actions of a decision-maker. The core components of the framework are: a state space $\mathcal{S}$, an action space $\mathcal{A}$, a transition mechanism $p_{t}:\mathcal{S}\times\mathcal{A}\times\mathcal{S}\to[0,1]$, a decision-making policy $\pi_{t}:\mathcal{S}\times\mathcal{A}\to[0,1]$, and an immediate cost function $C_{t}:\mathcal{S}\times\mathcal{A}\to\mathbb{R}^{+}$. The transition mechanism specifies the probabilities of moving from one state to another given a particular action, the policy is a probability distribution over the actions indicating how likely each action is to be chosen in a given state, and the immediate cost function assigns a cost to each action in each state. An MDP generates a temporally ordered sequence of tuples $\{(S_{t},A_{t},S_{t+1}):t\in\mathcal{T}\}$, modeling situations where the goal is to minimize cumulative costs over time.

The holding cost at an entry point \(i\) is computed using \(h_{i,t}:\mathbb{R}\to\mathbb{R}^{+}\), a convex and non-decreasing function of the state variable \(S_{i,t}\). At an exit point \(j\), the holding and backorder cost, denoted by \(h_{j,t}:\mathbb{R}\to\mathbb{R}^{+}\), is similarly a convex and non-decreasing function of \(S_{j,t}\). \cite{boujemaa2022sddp} considered a linear inventory cost function at entry points, defined as \(h_{i}(S_{i,t}) = \mathit{CW}_{i}S_{i,t}\), where \(\mathit{CW}_{i} \ge 0\) is a constant unit holding cost coefficient. For the inventory and backorder costs at exit points, they used a time-homogeneous cost function \(h_{j}(S_{j,t}) = \mathit{CD}_{j}S_{j,t}^{+} + \mathit{CB}_{j}S_{j,t}^{-}\), where \(\mathit{CD}_{j}>0\) and \(\mathit{CB}_{j}>0\) are constant unit inventory and backorder cost coefficients, respectively, with \(S_{j,t}^{+} = \max(S_{j,t},0)\) and \(S_{j,t}^{-} \doteq -\min(S_{j,t},0)\). Eventually, we define the total holding cost component as:
\begin{equation}
    h(S_{t}) = \sum_{i\in\mathcal{I}}h_{i}(S_{i,t})+\sum_{j\in\mathcal{J}}h_{j}(S_{j,t})
\end{equation}
We now define the immediate cost function that combines the total holding cost with the cost of executing a volume $A_{t}$ of drayage moves to external carriers under a portfolio contract. Let \(\check{A}_{t}\) denote the allocation vector \((A_{i,j,t}^{k} : (i,j) \in \mathcal{L}(k), k \in \mathcal{K})\), which represents the moves across the lanes by each carrier at time \(t\). Here, \(\mathcal{K} = \{1, \dots, n\}\) is the index set enumerating the elements of \(\mathit{CS} \cup \mathit{CO}\), where \(\mathit{CS}\) is the set of strategic carriers and \(\mathit{CO}\) is that of spot carriers. The symbol $\overline{S}_{j}$ indicates the maximum capacity for containers at each exit hub $j$. The immediate cost function is defined as follows:
\begin{align}
    C_{t}(S_{t},A_{t},Z_{t}) = h(S_{t})&+\min_{\check{A}_{t}}\sum_{k \in \mathcal{K}}\sum_{(i,j)\in\mathcal{L}(k)}W_{i,j,t}^{k}A_{i,j,t}^{k}\label{eq:imcostfst}\\
    \text{subject to}\hspace{5pt} & \sum_{k \in \mathcal{K}}\sum_{(i,j)\in\mathcal{L}(k)}A_{i,j,t}^{k} = A_{t} & \label{eq:imcostvol}\\
     & \sum_{(i,j)\in\mathcal{L}(k)}A_{i,j,t}^{k}\le x_{t}^{k} & \forall\,k\in\mathcal{K}  \label{eq:imcostccap}\\
     & \sum_{k \in \mathcal{K}}\sum_{j\in\mathcal{J}}A_{i,j,t}^{k} \le S_{i,t} + Q_{i,t} & \forall\,i\in\mathcal{I} \label{eq:imcostfromi}\\
     & \sum_{k \in \mathcal{K}}\sum_{i\in\mathcal{I}}A_{i,j,t}^{k} \le \overline{S}_{j} - S_{j,t} & \forall\,j\in\mathcal{J} \label{eq:imcosttoj}\\
     & A_{i,j,t}^{k} \ge 0 &\hspace{-20pt} \forall\,(i,j)\in\mathcal{L}(k),k\in\mathcal{K}\label{eq:imcostlst}
\end{align}
% Describe the constraints
The objective of the immediate cost function is to minimize the total transportation cost at time \( t \) by optimizing the allocation of volume across all lanes served by each carrier \( k \). The constraints ensure the following: Constraint~\eqref{eq:imcostvol} enforces that the total allocated volume matches the decision variable \( A_t \); Constraint~\eqref{eq:imcostccap} limits the total volume assigned to each carrier to its strategic capacity; Constraint~\eqref{eq:imcostfromi} guarantees that shipments from each entry location \( i \) do not exceed the available stock plus the inflow volume at time \( t \); Constraint~\eqref{eq:imcosttoj} restricts the volume received at each exit location \( j \) to its remaining capacity; finally, Constraint~\eqref{eq:imcostlst} ensures non-negativity for all allocation variables.

Fix an instance $\phi_{\tau}$ of the uncertainties that span the entire planning horizon; this instance is typically referred to as a scenario, representing a specific draw from the distribution of uncertain variables. We define the value $V_{1}^{\pi}(S_{1})$ of the random initial state $S_{1}$ under policy $\pi = (\pi_{t}:t\in\mathcal{T})$ as the expected total negative cost of operations:
\begin{equation}
    V_{1}^{\pi}(S_{1}) = -E_{\pi}\left[\sum_{t=1}^{\tau}C_{t}(S_{t},A_{t})\,\bigg|\,S_{1}\right]
\end{equation}
Note that we omit the scenario in the above notation, such that $V_{1}^{\pi}(S_{1}) = V_{1}^{\pi}(S_{1}, \phi_{\tau})$ and $C_{t}(S_{t}, A_{t}) = C_{t}(S_{t}, A_{t}, z_{t})$ for all $t \in \mathcal{T}$. Moreover, $E_{\pi}(\,\cdot\,)$ is an integration with respect to the probability measure:
\begin{equation}
    \prod_{t=1}^{\tau}p_{t}(s_{t+1} \mid s_{t},a_{t})\,\pi_{t}(a_{t} \mid s_{t})
\end{equation}
To obtain the optimal value of a cycle of operations, denoted by $V_{1}^{*}(S_{1})$, we need to execute the optimal policy $\pi^{*}$. To determine $\pi^{*}$, we propose the use of dynamic programming.

\subsubsection{Dynamic programming}\label{sec:dp}

We identify the optimal deterministic policy $\pi^{*}$ from a set of policies $\Pi$, using the backward induction method of dynamic programming:
\begin{align}\label{eq:dp}
    &V_{\tau+1}(s_{\tau+1})=-\alpha^{\top} s_{\tau+1}\\[5pt]
    &V_{t}(s_{t})=\max_{a_{t}\in\mathcal{A}(s_{t})}\left\{-C_{t}(s_{t},a_{t})+E_{S_{t+1}}V_{t+1}(s_{t},a_{t})\right\}\text{ for }t=\tau,\dots,1\label{eq:dpiter}
\end{align}

We assume the terminal value has a linear form with slope $\alpha \ge 0$ and associate the expected next state value function $E_{S_{t+1}}V_{t+1}:\mathcal{S}\times\mathcal{A}\to\mathbb{R}$ with the mapping $(s_{t},a_{t}) \mapsto \int_{\mathcal{S}} V_{t+1}(s_{t+1})\,p_{t}(ds_{t+1} \mid s_{t},a_{t})$ for all $t\in\mathcal{T}$. We note that the algorithm operates sequentially from time stage $\tau$ down to $1$. At each stage, it explores all possible states $s_t$ within the state space $\mathcal{S}$. For each state $s_{t}$, the algorithm computes the value of each action in the feasible action set, $\mathcal{A}(s_{t})$, and identifies the action $a_t^*$ that maximizes the value. This optimal action is then recorded in the policy function $\pi^{*}$, which is set to $1$ for the pair $(s_{t}^{\,}, a_{t}^{*})$ that achieves the maximum value. At the end of the algorithm, we have computed the optimal value $V_{1}^*(s_{1}^{\,})$ for any given starting state $s_{1}$.

\subsubsection{Computational complexity analysis}\label{sec:complexity}

The computational complexity of solving a $\tau$-horizon MDP with finite and discrete state and action spaces, where the cardinalities of these spaces are $|\mathcal{S}|$ and $|\mathcal{A}|$ respectively, to optimality using the backward induction algorithm is upper bounded by a finite multiple of $\tau|\mathcal{A}||\mathcal{S}|^2$, as noted by \citep{Rust1997}. The derivation of this expression can be understood through the steps of the backward induction algorithm: for each time period $t$, the algorithm iterates through all states. For each state $s_t$, it then considers every possible action from $\mathcal{A}(s_t)$, which is a subset of the action set $\mathcal{A}$. Given $\tau$ time periods, $|\mathcal{S}|$ possible states, and up to $|\mathcal{A}|$ possible actions per state, the algorithm conducts $\tau |\mathcal{A}| |\mathcal{S}|$ evaluations of immediate and future rewards, as reflected in the summation inside the $\max(\,\cdot\,)$ operator in Equation~\eqref{eq:dpiter}. Additionally, the computation of the expected future value for a given state-action pair involves summing over all possible next states $s_{t+1}$:
\begin{equation}
    E_{S_{t+1}}V_{t+1}(s_{t},a_{t}) = \sum_{s_{t+1}\in\mathcal{S}} V_{t+1}(s_{t+1})\,p_{t}(s_{t+1} \mid s_{t},a_{t})
\end{equation}
This sum involves $|\mathcal{S}|$ evaluations for each state-action pair, leading to the squared term $|\mathcal{S}|^2$ in the complexity expression. Each state is evaluated against every possible next state, thus justifying the quadratic dependence on the number of states in the computational complexity.

In our problem, the transition to a new state is governed by the mechanism of uncertain variables, as seen in Equation~\eqref{eq:draydyn}. We define the scenario space as the sample space of the probability space that models these uncertainties. Within this framework, uncertainties are termed exogenous state variables. These variables are not influenced by the decision-maker's actions but vary randomly in accordance with laws dictated by external forces. Understanding the cardinality of the scenario space is essential for accurately determining the computational complexity associated with our model.

The exogenous state of the system is represented by the random variable \(\Phi_t\), which is defined on the filtered probability space \((\Omega, \mathcal{F}, \{ \mathcal{F}_{t} : t \ge 0 \}, P_{t})\). If the exogenous state \(\Phi_t\) has a stationary distribution, the probability measure governing \(\Phi_t\) remains invariant under time shifts. Specifically, there exists a finite \(\Delta t\) such that for all \(t \geq 1\), the marginal distributions of \(\Phi_{\Delta t}\) and \(\Phi_{t+\Delta t}\) are identical. Under stationarity, the time-dependent probability measures \( P_{t} \) become equivalent across time shifts. However, the filtration \(\mathcal{F}_t\), representing the evolution of information up to time \(t\), may differ over time.

We assume that there exists a finite time window \(\Delta t > 0\) such that the distribution of the exogenous state becomes stationary for all \(t \geq \Delta t\). This implies that the probability measures governing the exogenous state, \(P_t\), remain invariant for \(t \geq \Delta t\). Consequently, the expected future value can be expressed as:
\begin{equation*}
\begin{aligned}
    E_{S_{t+1}}V_{t+1}(s_{t},a_{t}) &= \sum_{s_{t+1}\in\mathcal{S}} V_{t+1}(s_{t+1})\,p_{\Delta t}(s_{t+1} \mid s_{t},a_{t})\\
    &= \sum_{z_{t}\in\mathcal{Z}}\sum_{s_{t+1}\in\mathcal{S}} V_{t+1}(s_{t+1})\,\delta(s_{t+1}-f_{S,t}(s_{t},a_{t},z_{t}))\,p_{\Delta t}(z_{t})\\
    &= \sum_{z_{t}\in\mathcal{Z}} V_{t+1}(f_{S,t}(s_{t},a_{t},z_{t}))\,p_{\Delta t}(z_{t})
\end{aligned}
\end{equation*}
Deriving the expected future value requires evaluating the value function across the support $\mathcal{Z}$ of the uncertainty. The computational complexity of identifying an optimal volume allocation policy for the drayage problem described in this work is upper bounded by $\tau|\mathcal{A}||\mathcal{S}||\mathcal{Z}|$. In the remainder of this section, we investigate the individual components of the expression in greater detail.

% An exhaustive approach to the backward induction algorithm calls for looping over all possible $|\mathcal{Z}^\tau|$ scenarios. This approach captures uncertainties that, though exogenous to the system, become integral to the state variables through their influence on transition probabilities.

\subsubsection*{Cardinality of the endogenous state space}

Suppose $S_{i,t}\in\{0,1,\dots,\overline{S}_{i}\}$ for all $t\in\mathcal{T}$ and $i\in\mathcal{I}$, and $S_{j,t}\in\{-\underline{S}_{j},\dots,\overline{S}_{j}\}$ for all $t\in\mathcal{T}$ and $j\in\mathcal{J}$. We can easily compute the cardinality of the state space $\mathcal{S}$ with the following formula:
\begin{equation}
    |\mathcal{S}|=\prod_{i\in\mathcal{I}}(\overline{S}_{i}+1)\prod_{j\in\mathcal{J}}(\underline{S}_{j}+\overline{S}_{j}+1)
\end{equation}
Assuming that $\overline{S}_{i}=\overline{S}_{\mathcal{I}}$ for all $i\in\mathcal{I}$ and $(\underline{S}_{j},\overline{S}_{j})=(\underline{S}_{\mathcal{J}},\overline{S}_{\mathcal{J}})$ for all $j\in\mathcal{J}$, we obtain the simplified version $|\mathcal{S}|=(\overline{S}_{\mathcal{I}}+1)^{\mathcal{|I|}}(\underline{S}_{\mathcal{J}}+\overline{S}_{\mathcal{J}}+1)^\mathcal{|J|}$. This indicates that for fixed values of $\overline{S}_{\mathcal{I}}$, $\underline{S}_{\mathcal{J}}$, and $\overline{S}_{\mathcal{J}}$ the cardinality of the state space grows exponentially with the number of origin and destination points. For example, if $\overline{S}_{\mathcal{I}}=\overline{S}_{\mathcal{J}}=9$, $\underline{S}_{\mathcal{J}}=10$, and $|\mathcal{I}|=4$ and $|\mathcal{J}|=2$, we get $|\mathcal{S}|=4 \cdot 10^{6}$, whereas if, instead, $|\mathcal{I}|=6$ and $|\mathcal{J}|=3$, then $|\mathcal{S}|=8 \cdot 10^{9}$. Considering that the dynamic programming algorithm must loop through all states at each time stage, it is evident that the optimization of operations is computationally intensive even for small instances of the problem.

\subsubsection*{Cardinality of the scenario space support}

The exogenous state vector $\Phi_{t}$ remains impervious to the internal dynamics of our drayage system but serves as a stimulus, prompting the intervention of an internal agency to control the behavior of the system. In our problem, we define $\Phi_{t} \doteq \{\Phi_{t-1}, (Q_{i,t}: i \in \mathcal{I}), (D_{j,t}: j \in \mathcal{J}), (W_{i,j,t}^{k}: (i,j) \in \mathcal{L}(k), k \in \mathit{CO})\}$ as the exogenous state vector for $t\in\mathcal{T}$ and we let $\Phi_{0}\doteq\emptyset$. Suppose that for each $t\in\mathcal{T}$, $W_{i,j,t}^{k}$ takes a value from the set $\mathcal{W}$ for all $(i,j)\in\mathcal{L}(k)$ and $k\in\mathit{CO}$. Furthermore, for all $i\in\mathcal{I}$, $Q_{i,t}$ can take any value in $\mathcal{Q}_{i}$ and, similarly, $D_{j,t}$ takes values in $\mathcal{D}_{j}$ for all $j\in\mathcal{J}$. Hence, at any time $t\in\mathcal{T}$, the exogenous state variable can be any of the following number of instances:
\begin{equation}\label{eq:exogbound}
\begin{aligned}
    |\mathcal{Z}^{t}|&=\left(\prod_{i\in\mathcal{I}}|\mathcal{Q}_{i}|\prod_{j\in\mathcal{J}}|\mathcal{D}_{j}|\prod_{k\in\mathit{CO}}|\mathcal{W}|^{|\mathcal{L}(k)|}\right)^{t} \\
    &\le\left(|\mathcal{Q}|^{|\mathcal{I}|}|\mathcal{D}|^{|\mathcal{J}|}|\mathcal{W}|^{|\mathit{CO}||\mathcal{I}\times\mathcal{J}|}\right)^{t}
\end{aligned}
\end{equation}
To derive the above expression, we define $\mathcal{Q}\doteq\bigcup_{i\in\mathcal{I}}\mathcal{Q}_{i}$, $\mathcal{D}\doteq\bigcup_{j\in\mathcal{J}}\mathcal{D}_{j}$, and observe that $\mathcal{L}(k)\subset\mathcal{I}\times\mathcal{J}$ for all $k\in\mathit{CO}$. Overall, the scale of the exogenous state space is determined by the product of the components on the right hand side of~\eqref{eq:exogbound}. For instance, if $|\mathcal{I}|=|\mathcal{J}|=|\mathit{CO}|=2$ and $|\mathcal{Q}|=|\mathcal{D}|=|\mathcal{W}|=10$, we obtain a massive $10^{12t}$ possible values of the exogenous state vector at each $t$. Consequently, an exhaustive search over the entire space may quickly become computationally infeasible.

\subsubsection*{Cardinality of the action space}

At each time stage $t \in \mathcal{T}$, the decision-maker intervenes by instructing the transfer of a designated quantity $A_{t}$ of containers from entry to exit points via truck. This decision is made in consideration of both the endogenously evolving state of the system and the uncertainty encapsulated in the exogenous state. The objective is to select an action that optimizes the immediate reward, characterized by minimizing holding and transportation expenses, while simultaneously maximizing the expected value associated with the subsequent state transition following the execution of the decision. By introducing the portfolio contract described by \cite{Martinez2005}, the assignment of drayage moves among available carrier sources is internally executed within the computation of immediate transportation costs at each $t\in\mathcal{T}$, effectively narrowing our action space from a massive set encompassing all potential assignments between carriers to a finite set of discrete quantities, denoted $\mathcal{A}$.

\subsubsection{Approximate dynamic programming}

Let us define the Bellman operator under perfect information as the function $B(V):\mathcal{S}\to\mathbb{R}$ that, for some discount factor $\gamma\in(0,1]$, is associated with the mapping:
\begin{equation}
    s_{t}\mapsto\max_{a_{t}\in\mathcal{A}(s_{t})}\left\{\int_{\mathcal{Z}} \big[ -C_{t}(s_{t},a_{t},z_{t})+\gamma V_{t+1}(f_{S,t}(s_{t},a_{t},z_{t})) \big]\,p_{t}(dz_{t})\right\}
\end{equation}
We denote this as \( B(V_{t+1})(s_{t}) \). By \textit{perfect information}, we refer to the case where the decision-maker has access to the transition mechanism underlying the evolution of the system. In our problem, this means knowledge of the nuisance parameter of the transition mechanism, which is the density function of the exogenous state vector.

Using the Bellman operator with $\gamma=1$, we can rewrite Equation~\eqref{eq:dpiter} in the following compact form:
\begin{equation}\label{eq:backdp2}
    V_{\tau-t}=B^{t}(V_{\tau})
\end{equation}
Evaluating the integral in the Bellman operator is computationally intensive because, in the absence of a closed-form expression, we must sum over all possible realizations $z_{t}$ of $Z_{t}$ with respect to the known probability density $p_{t}(z_{t})$. Equation~\ref{eq:exogbound} illustrates how rapidly exact computation can become infeasible in the discrete case.

Suppose, we draw $N$ examples $\omega_{1}, \dots, \omega_{N}$ from $\Omega$ and let $z_{t,\ell} \doteq Z_{t}(\omega_{\ell})$ for all $t = 1, \dots, \tau$ and $\ell = 1, \dots, N$. We can approximate the integral in the Bellman operator by:
\begin{equation}
\begin{aligned}
    \int_{\mathcal{Z}} & \big[ -C_{t}(s_{t},a_{t},z_{t})+\gamma V_{t+1}(f_{S,t}(s_{t},a_{t},z_{t})) \big] \, p_{t}(dz_{t}) \\
    \approx & \sum_{\ell=1}^{N}\big[ -C_{t}(s_{t},a_{t},z_{t,\ell})+\gamma V_{t+1}(f_{S,t}(s_{t},a_{t},z_{t,\ell})) \big] \, \frac{p_{t}(z_{t,\ell})}{\sum_{l=1}^{N}p_{t}(z_{t,l})}
\end{aligned}
\end{equation}
Furthermore, we define the approximate Bellman operator $\widetilde{B}_{N}(V):\mathcal{S}\to\mathbb{R}$ under perfect information, for all $t\in\mathcal{T}$, $s_{t}\in\mathcal{S}$, and a given scenario $\phi_{\tau}$ as:
\begin{equation}
\begin{aligned}
    \widetilde{B}_{N}&(V_{t+1})(s_{t}) \doteq \\
    & \max_{a\in\mathcal{A}_{t}} \left\{ \sum_{\ell=1}^{N}\big[ -C_{t}(s_{t},a_{t},z_{t,\ell})+\gamma V_{t+1}(f_{S,t}(s_{t},a_{t},z_{t,\ell})) \big] \, \widetilde{p}_{t}(z_{t,\ell}) \right\}
\end{aligned}
\end{equation}
Here, we let $\widetilde{p}_{t}(z_{t,\ell})\doteq p_{t}(z_{t,\ell})/\sum_{l=1}^{N}p_{t}(z_{t,l})$ to ensure that the probability weight for each sample is properly normalized. By using the approximate method, we reduce the number of evaluations to \( \tau |\mathcal{A}| |\mathcal{S}| N \), which can prove highly beneficial when \( N \ll |\mathcal{Z}| \).

\subsection{Capacity optimization}

Let $\mathcal{X}$ represent the set of feasible capacities that a shipper can plan, at time $t = 0$, before the outset of operations. An element $x$ of $\mathcal{X}$ is a vector of capacity arrangements defined by:
\begin{equation}\label{eq:capvec}
    x\doteq\left(x_{t}^{k}\ge0:k\in\mathcal{K},t\in\mathcal{T}\right)
\end{equation}
Here, \(x_t^k\) denotes the non-negative capacity with each carrier \(k \in \mathcal{K}\) across all time periods \(t \in \mathcal{T}\). We assume that \(\mathcal{X}\) is convex. At time period \(t=0\), there is a cost associated with reserving capacities \(x \in \mathcal{X}\), which can be represented as a convex function \(v(x)\). According to \cite{Martinez2005}, a common choice is a cost function \(v\) that is linear in the capacities:
\begin{equation}
    v(x)\doteq\sum_{t=1}^{\tau}\sum_{k=1}^{n}v_{t}^{k}x_{t}^{k}
\end{equation}
It follows from Theorem 2 and Corollary 1 \citep{Martinez2005} that the problem of optimal capacity planning, addressed at the beginning of the planning horizon, is a concave optimization problem. Consequently, there exists a unique solution for selecting the optimal capacities from an initial state \(s_1\) prior to the commencement of operations, which can be formulated as follows:
\begin{equation}\label{eq:strategic}
    \max\left\{V_{1}^{*}(s_{1},x)-v(x):x\in\mathcal{X}\right\}
\end{equation}
The vector of carrier capacities, denoted as $x$, comprises individual components corresponding to carrier indices $k = 1, \dots, n$ and time period indices $t = 1, \dots, \tau$. To enhance clarity, we can represent the capacity parameters in a tabular format:
\begin{equation*}
\begin{matrix}
    x_{1}^{1} & x_{2}^{1} & \cdots & x_{\tau}^{1}\\
    x_{1}^{2} & x_{2}^{2} & \cdots & x_{\tau}^{2}\\
    \vdots & \vdots & \ddots & \vdots \\
    x_{1}^{n} & x_{2}^{n} & \cdots & x_{\tau}^{n}
\end{matrix}
\end{equation*}
Evidently, there are $n\tau$ capacities that need to be optimized in the first stage to address the strategic problem of optimal capacity reservation prior to operations. Consequently, for large values of $n$ and $\tau$, the search for optimal capacities $x^{*}$ may involve searching a highly-dimensional feasible set $\mathcal{X}$. This poses a significant challenge due to the computationally expensive nature of $V_{1}^{*}(S_{1},x)$. 

\section{Numerical experiments}\label{sec:numex}

To establish a case study for experimental analysis, several key parameters must be defined. First, determine the number of storage facilities where containers are held during entry or exit from the drayage system, which defines the index sets $\mathcal{I}$ and $\mathcal{J}$. Next, specify the total number of decision periods, denoted by $\tau$, to define the index set $\mathcal{T}$ for the planning horizon. The set of strategic carriers and their bids are essential for calculating transportation and capacity reservation costs. Each bid includes a subset of all possible lanes $\mathcal{L}$.

To conduct the experiment, the portfolio contract must be fully defined by specifying the set of sources, $\mathcal{K}$, which includes both strategic and spot market sources. The reservation cost $v_{t}^{k}$ for a unit of capacity $x_{t}^{k}$ associated with each source $k \in \mathcal{K}$ for each time period $t \in \mathcal{T}$ must be provided. A unit of capacity and a unit of transportation are measured in Twenty-foot Equivalent Units (TEUs). This capacity is reserved at the start of the planning horizon, before operations begin. Additionally, the execution cost $w_{i,j,t}^{k}$ for transporting a TEU across each lane $(i,j)$ under each source $k$ at each time $t \in \mathcal{T}$ must also be specified. In the case of a strategic source $k \in \mathit{CS}$, the unit transportation costs $\left( w_{i,j,t}^{k}:(i,j)\in\mathcal{L}(k) \right)$ are fixed and known from the start. However, for spot sources $\mathit{CO}$, a distribution $P_{W}$ must be specified to generate unit transportation costs $\left( w_{i,j,t}^{k} : (i,j) \in \mathcal{L}, k \in \mathit{CO}, t \in \mathcal{T} \right)$.

Finally, it is necessary to specify the holding costs associated with the temporary storage of TEUs at the entry and exit locations, as well as the distribution of TEU inflows and outflows, denoted by $P_{Q,D}$. This distribution governs the flow of TEUs into and out of the drayage system throughout the planning horizon, allowing us to simulate instances of  $\{ (Q_{t}, D_{t}) : t \in \mathcal{T}\}$. It is important to highlight that we have implicitly assumed the independence of random spot transportation costs from the random inflows and outflows. This assumption is valid, provided we condition on auxiliary information that eliminates any correlations between these variables.

\subsection{Instance generation logic}

The dynamic programming approach to the multistage container transportation problem under uncertainty is computationally intensive, as detailed in Section~\ref{sec:complexity}. This complexity is largely due to the expansive nature of both the state and scenario spaces. Specifically, we demonstrated that the state space increases exponentially with the number of entry and exit locations and grows polynomially with the storage capacities at these locations. Moreover, the scenario space's size also expands exponentially with the number of entry and exit locations and the number of spot sources, while it grows polynomially with the range of possible arrival or departure volumes and spot rates.

Collecting all the necessary data can be challenging, so we provide a method for generating synthetic instances. We fix the index sets $\mathcal{I}$ and $\mathcal{J}$, defining the set of possible lanes as $\mathcal{L} = \mathcal{I} \times \mathcal{J}$. A predetermined number of bids are generated by randomly selecting lanes from $\mathcal{L}$ and assigning a winner for each bid from a fixed set of carrier indices. The carriers winning these bids form the set of strategic sources. The number of spot sources is also fixed. For each strategic source, the unit transportation cost vector across lanes within the source's scope $\left( w_{i,j}^{k} : (i,j) \in \mathcal{L}(k), k \in \mathit{CS} \right)$ is sampled by independently drawing values from a common normal distribution, truncated below by a minimum value. Unit transportation costs are assumed to be homogeneous, independent across sources and lanes, and constant over time. The same logic applies to the unit transportation costs associated with spot sources $\left( w_{i,j,t}^{k} : (i,j) \in \mathcal{L}, k \in \mathit{CO}, t \in \mathcal{T} \right)$, which are also assumed to be homogeneous, independent across sources, lanes, and time periods. Finally, the capacities of strategic and spot sources are assigned fixed values, which may vary over time. However, it is intuitive that spot sources typically have larger capacities relative to strategic sources.

\subsection{Illustrative example}

To demonstrate the solution to the drayage problem, we examine a 4-period example ($\tau=4$) with a single entry and exit location. The example involves one strategic and one spot source, both of which service the sole lane between these locations. We denote the entry location index set as $\mathcal{I} = \{1\}$ and the exit location index set as $\mathcal{J} = \{|\mathcal{I}|+1\} = \{2\}$. Consequently, the set of lanes is $\mathcal{L} = \{(1, 2)\}$, representing the only lane from entry point $1$ to exit point $2$. The strategic carrier set, $\mathit{CS} = \{1\}$, and the spot carrier set, $\mathit{CO} = \{|\mathit{CS}|+1\} = \{2\}$, each contain a single element. Therefore, the set of all sources is $\mathcal{K} = \{1, 2\}$, as there are exactly two sources in this example.

The execution cost for the strategic source is fixed at \( w_{1,2}^{1} = \$14.7 \) per TEU for each period \( t = 1, \dots, \tau \), while the spot source cost for each \( t \) follows a distribution given by \( 0.4\,\delta(w_{1,2,t}^{2} - 7) + 0.6\,\delta(w_{1,2,t}^{2} - 22) \), meaning the cost is \( \$7.0 \) per TEU with a probability of 0.4 and \( \$22.0 \) with a probability of 0.6. For all \( t = 1, \dots, \tau \), the inflow of containers \( Q_{1,t} \) follows a distribution described by \( 0.4\,\delta(q_{1,t}) + 0.3\,\delta(q_{1,t} - 4) + 0.3\,\delta(q_{1,t} - 8) \), while the outflow \( D_{2,t} \) is distributed as \( 0.25\,\delta(d_{2,t}) + 0.25\,\delta(d_{2,t} - 4) + 0.5\,\delta(d_{2,t} - 8) \). The support of \( W_{1,2,t}^{2} \) is \( \mathcal{W} = \{7, 22\} \) for all \( t \), and the supports of \( Q_{1,t} \) and \( D_{2,t} \) are \( \mathcal{Q} = \mathcal{D} = \{0, 4, 8\} \) for all \( t \). Consequently, the support of the exogenous state variable \( Z_{t} \) is \( |\mathcal{Z}| = 18 \) for each period \( t \), leading to \( |\mathcal{Z}|^{\tau} = 104,976 \) potential scenarios when \( \tau = 4 \). This highlights the significant complexity that arises as the planning horizon increases, even in small instances of the problem.

\subsubsection{Operational volume allocation policy optimization}

To address the problem exhaustively, it is necessary to determine the optimal policy for each of the 104,976 possible scenarios. For a particular scenario $\phi_{4}=(z_{1},\dots,z_{4})$, this is achieved by iterating backward from \( t = 4 \) to \( t = 1 \) for all \( s_{t} \in \mathcal{S} \) using the following procedure:
\begin{equation}\label{eq:dpscen}
    V_{t}(s_{t}) = \max_{a_{t} \in \mathcal{A}(s_{t})} \left\{ -C_{t}(s_{t}, a_{t}, z_{t}) + V_{t+1}(f_{S,t}(s_{t}, a_{t}, z_{t})) \right\}
\end{equation}
where the terminal value is defined as \( V_{5}(s_{5}) = -\alpha^{\top} (s_{1,5}^{\,}, s_{2,5}^{+}, s_{2,5}^{-}) \). Here \(\alpha = (15, 12, 24)\) represents the holding cost in \$/TEU, with $\mathit{CW}_{1}=15$, $\mathit{CD}_{2}=12$, $\mathit{CB}_{2}=24$, and \( s_{\bullet}^{+} = \max(s_{\bullet}^{\,}, 0) \) and \( s_{\bullet}^{-} = \min(s_{\bullet}^{\,}, 0) \). We set the storage limits as \(\overline{S}_1 = \overline{S}_2 = \underline{S}_2 = 10\), resulting in a total of \( |\mathcal{S}| = 11(10 + 10 + 1) = 231 \) states. Additionally, the maximum TEU volume that can be moved within the drayage system is fixed to 10, which implies that the action set \(\mathcal{A}(s_t)\) is a subset of \(\{0, 1, \dots, 10\}\) for all \(s_t \in \mathcal{S}\) and for each \(t = 1, \dots, 4\). 

To gain a deeper understanding of the characteristics of the value surface, we first select a capacity arrangement and define a specific scenario. We then solve this scenario to optimality and plot the value function at each time step \( t = 1, \dots, 4 \). The chosen capacity arrangement and scenario details are presented in Table~\ref{tab:testinstance}. The initial state of the drayage system is set to \( S_{0} = (0, 8) \), indicating that the stock level at the entry is 0, while there are 8 TEUs at the exit.

\begin{table}[ht]
    \caption{Description of the reserved capacities and selected scenario.}
    \label{tab:testinstance}
    \centering
    \begin{tabular}{l|lrrrrl}
        \toprule
        & $t$ & 1 & 2 & 3 & 4 &  \\
        \midrule
        \multirow{2}{*}{\shortstack{Imposed \\ capacities}} & $x_{t}^{1}$ & 4 & 3 & 2 & 4 & TEU \\
        & $x_{t}^{2}$ & 4 & 4 & 4 & 4 & TEU \\
        \midrule
        \multirow{3}{*}{\shortstack{Uncertainty \\ realization}} & $q_{1,t}$ & 8 & 8 & 0 & 0 & TEU \\
        & $d_{2,t}$ & 8 & 8 & 8 & 0 & TEU \\
        & $w_{1,2,t}^{2}$ & 7 & 22 & 7 & 22 & \$/TEU \\
        \bottomrule
    \end{tabular}
\end{table}

\begin{figure}[ht]
    \centering
    \includegraphics[width = 0.85 \textwidth]{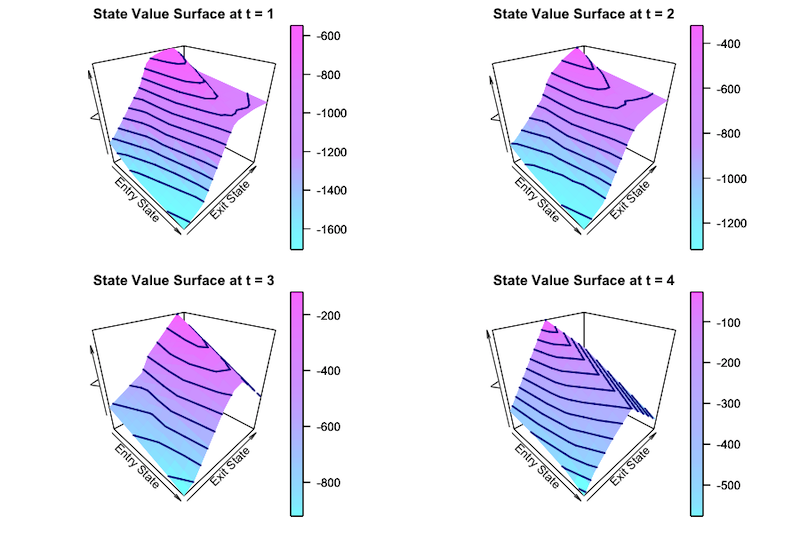}
    \caption{The value function for the chosen scenario is evaluated at each period of the planning horizon. The entry state ranges from $0$ to $10$, while the exit state varies between $-10$ and $10$, both defined as integer values. % Although the plot displays continuous surfaces for visual clarity, the value function itself is calculated strictly at these integer points. Contour lines are incorporated to emphasize the characteristics of each surface. 
    The color bar indicates the negative cumulative cost from the given period up to $t=4$.}
    \label{fig:i1x1Vsurf}
\end{figure}

Figure~\ref{fig:i1x1Vsurf} presents the value function for each state of the system, revealing its quasi-concave shape. Specifically, the function exhibits concavity at higher contour levels and convexity at lower ones. This pattern results from the linear holding cost function, which is designed with a steeper slope for negative exit states (backorders). At time points 1-3, the maximum value occurs with an entry state of 0 and exit states of 6 at time point 1 and 8 at time points 2 and 3. By time point 4, the maximum value is reached when both entry and exit states are 0, indicating a shift from maintaining stock at the exit location to no stock at either location.

\begin{figure}[ht]
    \centering
    \includegraphics[width = \textwidth]{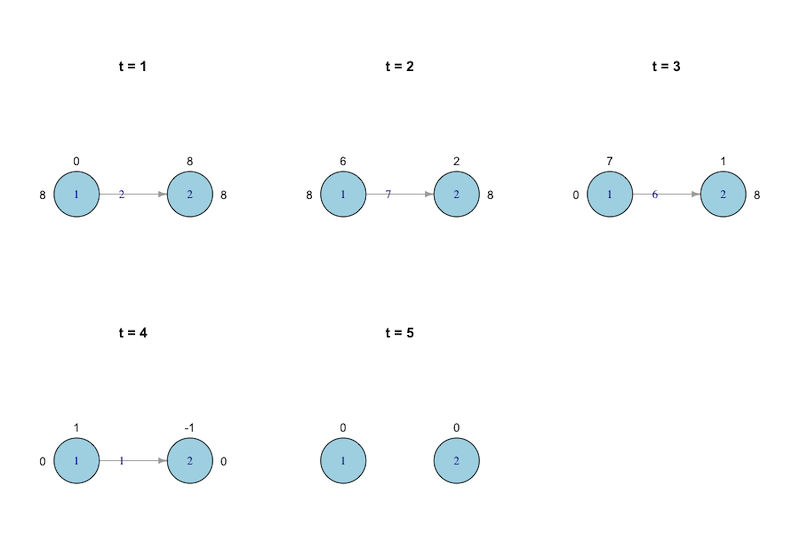}
    \caption{Evolution of the drayage system under the identified optimal volume allocation policy. The left node shows the entry location, and the right node indicates the exit. The top number on each node represents stock levels, while the left of the entry node shows incoming container volume, and the right of the exit node shows outflow. The middle number on the connecting arrow indicates the volume transported between the locations.}
    \label{fig:i1x1evolution}
\end{figure}

We demonstrate the evolution of the system under the optimal policy for the selected scenario, as shown in Table~\ref{tab:testinstance}, in Figure~\ref{fig:i1x1evolution}. TThe system's optimal value for this scenario is attained at the initial state \( s_{1}^{*} \doteq \argmax\{V_{1}^{*}(s_{1}, x) : s_{1} \in \mathcal{S}\} = (0,8) \), ensuring that the system begins with stock at the exit location. As time progresses, the optimal policy prioritizes depleting stock at the entry by transferring TEUs to the exit location. By the end of the planning horizon, the policy eliminates all remaining stock, aligning with the value function surfaces that show maximum value (or least cost) when both entry and exit states are minimized.

\subsubsection{Capacity strategy optimization}

We conclude this example by discussing the joint optimization of capacity $\{(x_{t}^{1},x_{t}^{2}): t = 1, \dots, 4\}$ and the total operational cost of the drayage system under the uncertainty realization detailed in Table~\ref{tab:testinstance}. The strategic source has a fixed reservation cost-per-TEU of $(v_{1}^{1}, v_{2}^{1}, v_{3}^{1}, v_{4}^{1}) = (\$8.52, \$4.46, \$4.25, \$9.70)$, with a uniform execution cost of $w_{1,2}^{1}=\$2.94$/TEU for all time periods $t=1,\dots,4$. No prior reservations are made for the spot carrier. The optimization is carried out by solving Problem~\eqref{eq:strategic}, assuming linear reservation costs. 

An iterative capacity search algorithm updates the capacity arrangement at each step by solving the dynamic program~\eqref{eq:dpscen}. However, since the goal is to determine the system’s optimal value at its initial state, a more computationally efficient approach is adopted. The problem is reformulated as a multistage linear program, following \citep{boujemaa2022sddp}. This relaxation enables the use of the L-BFGS-B algorithm via, for instance, the \texttt{optim()} function in \texttt{R} \citep{Rman}. L-BFGS-B belongs to the class of quasi-Newton methods, a subset of numerical optimization algorithms that approximate second-order information (the Hessian matrix) using only first-order gradient evaluations \citep{nocedal2006}. 

The objective value $V_{1}^{*}(x)-v(x)$ is evaluated at both the initial capacity arrangement $x^{0}$, shown in Table~\ref{tab:testinstance}, and the optimized capacity $x^{*}$, in Table~\ref{tab:optcap}, yielding values of \$557.2 and \$439.2, respectively. This reflects a 21.2\% reduction in total costs under the optimized capacity arrangement. Here, $V_{1}^{*}(x)$ represents the optimal value for a given capacity vector $x$ at the initial state $s_{1}$, maximizing the system's value at $t=1$, i.e., $V_{1}^{*}(x) = \max\{V_{1}^{*}(s_{1}, x) : s_{1} \in \mathcal{S}\}$. While we have demonstrated that $x^{*}$ yields a higher value than $x^{0}$, further investigation is necessary to conclusively establish that $x^{*}$ represents the optimal capacity strategy.

\begin{table}[ht]
    \caption{Optimal capacity arrangement for the selected scenario.}
    \label{tab:optcap}
    \centering
    \begin{tabular}{l|rrrr|l}
        \toprule
        $t$ & 1 & 2 & 3 & 4 & Unit \\
        \midrule
        $x_{t}^{1,*}$ & 0 & 8 & 0 & 0 & TEU \\
        $x_{t}^{2,*}$ & 4 & 4 & 8 & 4 & TEU \\
        \bottomrule
    \end{tabular}
\end{table}

In this example, a capacity arrangement consists of one strategic and one spot source over four time periods. With each source’s capacity in a given time period ranging from 0 to 10, the total number of possible capacity vectors is $(11^2)^4 = 11^8 = 214,358,881$. Computing the total cost (capacity reservation and operations) for each vector requires optimizing the value function, which is computationally expensive. To address this, we perform a Monte Carlo simulation by uniformly sampling 1,000,000 capacity vectors, calculating the total cost for each. Additionally, we evaluate another important performance measure, namely, the cost-per-TEU. Table~\ref{tab:mcsummary} shows that the solution $x^{*}$ minimizes both the total cost and cost-per-TEU, with the latter computed as $10.98$ \$/TEU. 

\begin{table}[ht]
  \centering
  \caption{Summary statistics for Cost-per-TEU and Total Cost based on $1,000,000$ Monte Carlo simulations of the capacity vector.}
  \label{tab:mcsummary}
  \begin{tabular}{l|rrrrrr}
   \toprule
   Statistic & Min & 1st Qtl & Median & Mean & 3rd Qtl & Max \\
   \midrule
   Cost-per-TEU & 10.98 & 13.28 & 14.25 & 14.38 & 15.34 & 21.81 \\
   Total Cost   & 439.2 & 527.7 & 566.2 & 579.6 & 612.6 & 1,671.5 \\
   \bottomrule
  \end{tabular}
\end{table}

\begin{figure}[ht]
    \centering
    \includegraphics[width = \textwidth]{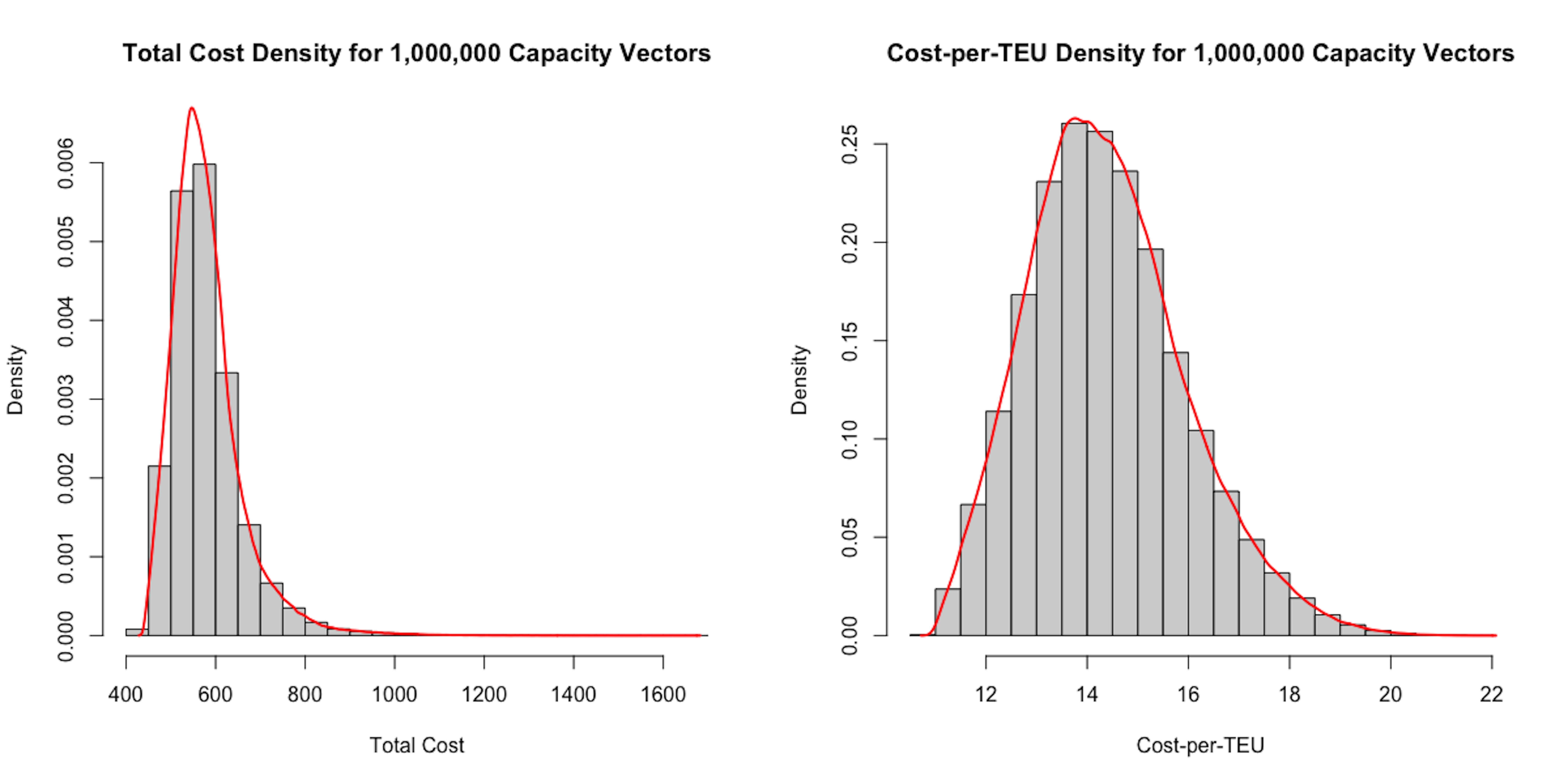}
    \caption{Density of Total Cost (left) and Cost-per-TEU (right) computed by optimizing the value function over $1,000,000$ random examples of the capacity arrangement.}
    \label{fig:dplots}
\end{figure}

Figure~\ref{fig:dplots} presents density plots of the simulated total costs and cost-per-TEU, illustrating a unimodal, slightly right-skewed distribution.The mode is located at considerably higher values than the minimum (leftmost) attained by $x^{*}$. The optimal strategy, reflected in Table~\ref{tab:optcap}, involves concentrating most capacity at the second and third time periods, primarily due to the high inflows and outflows during these intervals, as well as the presence of initial stock at the exit location. The absence of inflows or outflows in the final period minimizes the need for capacity adjustments in the later stages of the planning horizon. Furthermore, the strategy incorporates a buffer in the spot market, although capacity may occasionally remain unused. This underutilization is not problematic, as no reservation costs are incurred when planning capacity with spot carriers, thereby providing operational flexibility without additional financial burden.

\subsubsection{Approximate capacity optimization across scenarios}

Our analysis thus far has been confined to a specific scenario, allowing us to gain insights into the operational volume allocation policy and the capacity strategy under the given conditions. In this section, we broaden our approach by computing the capacity strategy that maximizes the value function across a subset of the scenario space. The optimal value of the initial state of the system $V_{1}^{*}(x)$, at some $x\in\mathcal{X}$, is now estimated by $\widehat{V}_{1,N}^{*}(x)$ using backward induction, expressed as \( V_{\tau - t} = \widetilde{B}_{N}^{t}(V_{\tau}) \), for $t = \tau, \dots, 1$, where \( \widetilde{B}_{N}(V_{\bullet}) \) is the approximate Bellman operator with discount factor \( \gamma = 1 \), applied to a sample of size \(N\), drawn from the set of all possible scenarios (on the order of 100,000).

We sample $N=1,000$ scenarios from $\mathcal{Z}^{4}$, assuming both stagewise independence and independence between inflow, outflow, and spot rates. Under these assumptions, the probability of each scenario is determined by the product of the probabilities of the individual time-point realizations, and each time-point realization's probability is the product of the probabilities of its component factors. While this study does not focus on the effects of correlations among uncertainties, incorporating such correlations would primarily influence the calculation of scenario probabilities, which could, in turn, affect the overall results.

\begin{figure}[ht]
    \centering
    \includegraphics[width = 0.75\textwidth]{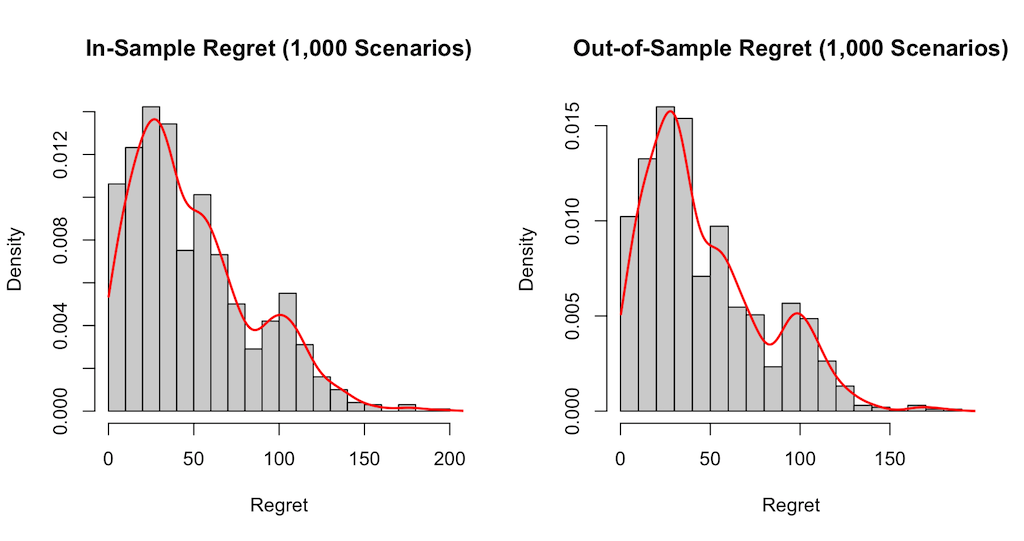}
    \caption{Regret density for the 1,000 scenarios used in optimizing the capacity strategy (left) and for 1,000 alternative scenarios (right).}
    \label{fig:regretdensity}
\end{figure}

For a given capacity strategy, represented by \(x \in \mathcal{X}\), we estimate the optimal value \(V_{1}^{*}(x)\) using the expected value functional over a sample of 1,000 scenarios. We then search for the capacity vector \(\widehat{x}_{N}^{\,*}\) that maximizes the value function approximation \(\widehat{V}_{1,N}^{*}(x)\). This is accomplished by applying a linear relaxation of the problem and utilizing the L-BFGS-B algorithm to find the solution. To evaluate the quality of the approximation \(\widehat{x}_{N}^{\,*}\), we define the regret for an individual scenario indexed by \(\ell\) as the difference between the optimal value and the value attained by the approximation. Specifically, the regret is given by:
\[
\text{Regret}_\ell = \left[ V_{1,\ell}^{*}(x_{\ell}^{*}) - v(x_{\ell}^{*}) \right] - \left[ V_{1,\ell}^{*}(\widehat{x}_{N}^{\,*}) - v(\widehat{x}_{N}^{\,*}) \right]\quad\text{for }\ell=1,\dots,N
\] % x.ev is 0-8 6-6 6-7 0-8
Here, \(V_{1,\ell}^{*}(x_{\ell}^{*})\) represents the optimal value for the given scenario \(\ell\) and its corresponding optimal strategy \(x_{\ell}^{*}\), while \(V_{1,\ell}^{*}(\widehat{x}_{N}^{\,*})\) is the value attained by the approximation \(\widehat{x}_{N}^{\,*}\) in that same scenario. The term \(v(x)\) is a reference value subtracted from both the optimal and approximate values to compute the regret. Figure~\ref{fig:regretdensity} illustrates the density of the regret computed over the same sample of scenarios used to obtain the value function approximation, as well as the density of the regret on scenarios outside of this sample. We observe the highest mode of the regret density at a value around \$30, with a second, lower mode at \$100. Notably, the scenarios corresponding to this second mode involve high inflow and outflow volumes, as well as elevated spot rates in most periods. In these cases, \(\widehat{x}_{N}^{\,*}\) did not plan capacity with strategic sources but instead relied on spot capacity, which was optimal in expectation but not for these particular instances.

\begin{figure}[ht]
    \centering
    \includegraphics[width = 0.5\textwidth]{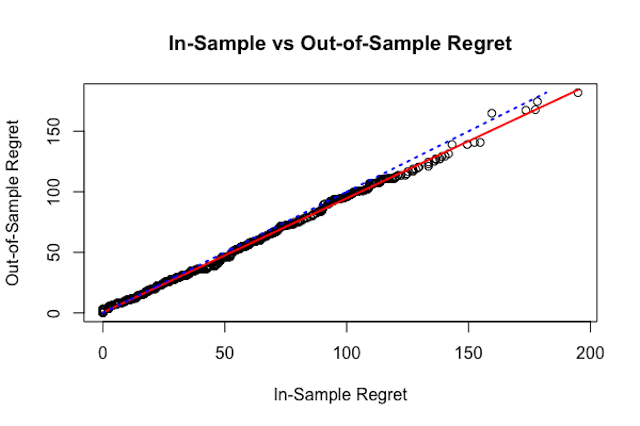}
    \caption{In-sample vs. out-of-sample regret plotted compared to the diagonal (dotted).}
    \label{fig:invsoutregret}
\end{figure}

Lastly, we investigate the extent to which the performance of the capacity vector \(\widehat{x}_{N}^{*}\) in maximizing the starting value of the system generalizes to new out-of-sample scenarios. To assess this, we compare the in-sample and out-of-sample regret by plotting them against each other in Figure~\ref{fig:invsoutregret}. The fact that the values align closely with the diagonal suggests that the regret in the new scenarios is very similar to that in the original sample, indicating strong generalization performance. It is worth noting that this result is achieved despite using only 1,000 scenarios to approximate the value function, out of a potential 100,000. This demonstrates that \(\widehat{x}_{N}^{*}\) effectively generalizes beyond the training sample, suggesting robustness in the optimization approach.

\subsection{Scalability considerations}
% Pave the way to SDDP application

The computational complexity of the dynamic programming algorithm, as discussed in Section~\ref{sec:optops}, is upper bounded by \( \tau|\mathcal{S}||\mathcal{A}||\mathcal{Z}| \) in the exact case, or by \( \tau|\mathcal{S}||\mathcal{A}|N \) in the approximate case. Let \( c \in (0, \infty) \) denote a constant such that the precise number of computations required by the algorithm is given by \( c\tau|\mathcal{S}||\mathcal{A}||\mathcal{Z}| \). The constant \( c \) primarily depends on the computational effort needed to determine the next state and immediate cost for a given state-action pair. To improve efficiency, all possible transitions and immediate rewards—i.e., tuples of the form (next state, immediate cost, present state, action)—can be precomputed and stored. This allows the dynamic programming algorithm to access transitions with minimal latency, reducing \( c \) to the time required for hashing the correct transition. The set of all such transitions fully defines the transition dynamics, commonly referred to as the ``model.''

Efficient model computation for the drayage problem is computationally intensive due to the need to solve a linear program \( \tau|\mathcal{S}||\mathcal{A}||\mathcal{W}||\mathcal{Q}| \) times to calculate immediate costs across all time periods, states, actions, spot rates, and inflow levels. As a result, even small instances require significant computational resources to solve optimally using dynamic programming. This challenge, known as the ``curse of dimensionality,'' is a fundamental limitation of the MDP framework. We propose using multistage linear programming relaxations when the goal is to determine the optimal value rather than the optimal policy, which is particularly advantageous when optimizing capacity as part of a joint optimization problem. Finally, it is important to note that as the number of carriers/sources $(n)$ and time periods $(\tau)$ increases, the dimensionality of the capacity vector grows significantly, making the search for an optimal capacity strategy increasingly challenging.

A straightforward strategy to reduce the dimensionality from $n\tau$ to $3n$ involves employing a quadratic parameterization for a sources' capacity across the operational horizon. Specifically, for all sources $k = 1, \dots, n$, we express it as:
\begin{equation}
x_{t}^{k} = \beta_{0}^{k} + \beta_{1}^{k}t + \beta_{2}^{k}t^{2}
\end{equation}
This approach evidently reduces the parameters to $3n$; however, it comes at the cost of foregoing the capability to capture more flexible temporal relationships. % cite Martinez-de-Albeniz

\section{Discussion}

This study introduced a joint optimization model integrating strategic capacity planning with operational volume allocation in drayage procurement. The primary objective was to bridge the gap between strategic planning and operational decision-making. Employing an MDP framework, we developed a dynamic programming-based approach to determine the optimal assignment of transportation volumes to carriers in a drayage environment characterized by uncertain container flows and spot market rates. The dynamic programming solution was then used to evaluate different capacity arrangements in search of the optimal one using the L-BFGS-B algorithm. 

While the computational experiments were conducted on a controlled example, the methodology was designed to be applicable to practical problems. A key contribution of this research is the explicit differentiation between carrier selection and capacity planning, which facilitates the effective integration of strategic and operational decisions. The computational complexity associated with this approach could be prohibitive, particularly due to the exponential growth of the state and scenario spaces. However, the adoption of approximate dynamic programming can significantly reduce computational demands, demonstrating its potential for larger-scale practical implementations. 

The experimental results provided insights into balancing strategic and spot carriers to handle varying container flows effectively. The optimal capacity arrangement involved reserving resources across both carrier types. To assess the broader applicability of the approach, we estimated the optimal capacity plan using 1,000 sampled scenarios and evaluated its performance against out-of-sample instances. The observed low regret in both in-sample and out-of-sample evaluations suggests that the proposed methodology generalizes well beyond the studied instance. Future research should further validate these findings in settings with industry-scale data and operational constraints.

Despite these contributions, the model has certain limitations. It assumes that uncertainties are independent at each time point and treats inflows and outflows as separate variables, which may not fully capture the interdependencies present in real-world scenarios. In practice, these variables may exhibit correlations, and future research should aim to extend the model to account for these relationships. Additionally, further refinement of the approximate dynamic programming approach could enhance the accuracy of policy approximations in larger and more complex scenarios.

Future work could develop more scalable and efficient approximation algorithms beyond the sample average approximation, accompanied by a more detailed modeling of the uncertain parameters. It could also explore the integration of real-time data streams into the decision-making process, enabling companies to continuously update their capacity plans and volume allocation policies based on current market and operational information. Additionally, incorporating machine learning techniques could enhance the modeling of uncertainties, leading to more accurate and adaptive policies.

In conclusion, this study establishes a foundation for enhancing strategic decisions in drayage procurement. By bridging the gap between strategy and operations, the proposed approach paves the way for more flexible, cost-efficient, and resilient logistics systems capable of adapting to the dynamic nature of global trade. As part of this effort, we are developing an \texttt{Rcpp} package to facilitate the implementation of the proposed models. The source code for the methods employed in this study is publicly available at \url{https://github.com/georgios-vassos1/TLPR}.

% \appendix

% \section{Sample Appendix Section}
% \label{sec:sample:appendix}
% Lorem ipsum dolor sit amet, consectetur adipiscing elit, sed do eiusmod tempor section \ref{sec:sample1} incididunt ut labore et dolore magna aliqua. Ut enim ad minim veniam, quis nostrud exercitation ullamco laboris nisi ut aliquip ex ea commodo consequat. Duis aute irure dolor in reprehenderit in voluptate velit esse cillum dolore eu fugiat nulla pariatur. Excepteur sint occaecat cupidatat non proident, sunt in culpa qui officia deserunt mollit anim id est laborum.

\bibliographystyle{elsarticle-num}
\bibliography{inland-refs}





% \end{thebibliography}
\end{document}